\documentclass[11pt]{article}
\usepackage{amssymb,amsmath,bbm,latexsym,amscd}
\newtheorem{theorem}[equation]{Theorem}
\newtheorem{corollary}[equation]{Corollary}

\newtheorem{remark}[equation]{Remark}
\numberwithin{equation}{section}

\newcommand{\fg}{\mathfrak{g}}

\newcommand{\ot}{\otimes}

\newcommand{\fh}{{\mathfrak h}}

\newcommand{\ga}{\alpha}

\newcommand{\proof}{{\bf Proof\ \ }}
\newcommand{\qed}{\hfill $\Box$}

\newcommand{\mf}{\mathfrak}

\title{Whittaker vectors in singular Whittaker modules}

\author{Karthik Dulam, Hrishikesh Ghate, Michael Lau, and Suyash Pathak}

\date{}
\begin{document}
\maketitle

\begin{small}
  \noindent {\bf Abstract:} Let $\fg$ be a complex semisimple Lie algebra with Borel subalgebra $\mf{b}$ and corresponding nilradical $\mf{n}$.  We show that singular Whittaker modules $M$ are simple if and only if the space $\hbox{Wh}\,M$ of Whittaker vectors is $1$-dimensional.  For arbitrary locally $\mf{n}$-finite $\fg$-modules $V$, an immediate corollary is that the dimension of $\hbox{Wh}\,V$ is bounded by the composition length of $V$.

  %For any finite-dimensional simple Lie algebra $\fg$ and finitely generated commutative associative algebra $S$, we give a complete classification of the simple weight modules of $\fg\ot S$ with bounded weight multiplicities. 

\bigskip


%weight modules, current algebras, admissible representations, infinite-dimensional Lie algebras 

\bigskip

\end{small}

\section{Introduction}

Let $\fg$ be a complex semisimple Lie algebra with Cartan subalgebra $\fh$, root system $\Phi$, base $\Delta$ of simple roots, and corresponding set $\Phi_+$ of positive roots.  The nilradical $\mathfrak{n}=\bigoplus_{\ga\in\Phi_+}\fg_\ga$ of the standard Borel acts trivially on the cyclic vectors of highest weight representations.  More generally, we consider Whittaker modules $V$ of type $\psi$, cyclic modules generated by a vector $v\in V$ on which $\mathfrak{n}$ acts by a Whittaker character $\psi:\ \mathfrak{n}\rightarrow \mathbb{C}$.  When $\psi(e_\ga)\neq 0$ for all $\ga\in\Delta$, the module $V$ and the character $\psi$ are said to be {\em non-singular}.

Whittaker modules were introduced by Kostant to address questions about primitive ideals, representations of semisimple Lie groups, and Toda integrable systems \cite{kostant78,kostant79}.  They later played a prominent role in Block's classification \cite{block80} of all irreducible representations of $\mathfrak{sl}_2$.  Recent interest has grown significantly due to equivalences between categories of generalized Whittaker modules and modules for finite $W$-algebras \cite{skryabin02}.  Under this correspondence, a generalized Whittaker module $M$ corresponds to its space $\hbox{Wh}\,M$ of Whittaker vectors, viewed as a module over a finite $W$-algebra.

The purpose of this short paper is to show that a Whittaker module $M$ is simple if and only if its space of Whittaker vectors is $1$-dimensional, even when $M$ is singular.  This recovers a result of Kostant in the non-singular case \cite[Theorem 3.6.1]{kostant78}, and may be somewhat surprising, as the associated $W$-algebras are generally non-commutative.

For arbitrary locally $\mathfrak{n}$-finite $\fg$-modules $M$, an immediate corollary is that the dimension of $\hbox{Wh}\,M$ is bounded by the composition length of $M$.  In particular, the dimension of the space of Whittaker vectors in a Whittaker module with central character $\chi$ and Whittaker character $\psi$ is always bounded by the length of the universal module with these properties.  The lengths of these universal modules are known, thanks to work by Mili\v{c}i\'c and Soergel \cite{milicic-soergel97} using Kazhdan-Lusztig theory.

\bigskip

\noindent
{\bf Acknowledgements.} This paper is the result of a MITACS Globalink undergraduate research internship held by Dulam, Ghate, and Pathak in Summer 2023 under the direction of Lau at Universit\'e Laval.  Financial and logistical support of MITACS and Universit\'e Laval are gratefully acknowledged.  The authors also thank Charlotte Lavoie-Bel for many useful discussions throughout the summer.

\section{Simple Whittaker modules}

We recall McDowell's classification of simple singular Whittaker modules \cite{mcdowell85}.  Let $\fg$, $\fh$, $\Phi$, $\Phi_+$, and $\Delta$ be as in Section 1, and let $\mathfrak{n}$ be the nilradical of the standard Borel $\mathfrak{b}$ relative to $(\fg,\fh,\Delta)$.  A Lie algebra homomorphism $\psi:\ \mathfrak{n}\rightarrow \mathbb{C}$ is called a {\em Whittaker character}; the character $\psi$ is {\em non-singular} if it is nonzero on the root spaces $\fg_\ga$ for all simple roots $\ga\in\Delta$.  An element $v$ of a $\fg$-module $V$ is said to be a {\em Whittaker vector of type} $\psi$ if $xv=\psi(x)v$ for all $x\in\mathfrak{n}$.  The space of Whittaker vectors in $V$ will be denoted by $\hbox{Wh}\,V$.  Cyclic modules generated by Whittaker vectors of type $\psi$ are called {\em Whittaker modules of type} $\psi$ and are {\em singular} or {\em non-singular}, depending on whether $\psi$ is singular or non-singular.  Each Whittaker module admits a unique Whittaker character, so the notions of singular and non-singular are well defined.

Let $\psi$ be a Whittaker character, and let $\Delta_\psi$ be the set of simple roots $\ga$ for which $\psi(\fg_\ga)\neq 0$. We write $\mathfrak{l}=\mathfrak{l}_\psi$ for the reductive subalgebra generated by the sum of $\fh$ and the root spaces $\fg_\ga$ with $\pm\ga\in \Delta_\psi$.  The Lie algebra $\mathfrak{l}$ has a triangular decomposition $\mathfrak{l}=\mathfrak{l}_-\oplus\fh\oplus\mathfrak{l}_+$, where $\mathfrak{l}_+$ and $\mathfrak{l}_-$ are the intersections of $\mathfrak{l}$ with the positive and negative root spaces of $\fg$, respectively.  The centre $\mathfrak{z}$ of $\mathfrak{l}$ lies in $\fh$, and the corresponding parabolic subalgebra $\mathfrak{l}+\mathfrak{n}$ will be denoted by $\mathfrak{p}$.  The subalgebra $\mathfrak{p}$ has a unique $\operatorname{ad}\fh$-stable decomposition $\mathfrak{p}=\mathfrak{l}\oplus\mathfrak{m}$, and $\fg=\overline{\mathfrak{m}}\oplus\mathfrak{l}\oplus\mathfrak{m}$, where $\overline{\mathfrak{m}}$ is the Lie subalgebra spanned by the root spaces $\fg_{-\ga}$ for which $\fg_\ga$ is contained in $\mathfrak{m}$.  For any Lie algebra $\mathfrak{a}$, let  $Z(\mathfrak{a})$ be the centre of its enveloping algebra $U(\mathfrak{a})$.  For an arbitrary central character $\Omega:\,Z(\mathfrak{l})\rightarrow\mathbb{C}$, consider the induced $\mathfrak{l}$-module
$$Y_{\Omega,\psi}=U(\mathfrak{l})\ot_{Z(\mathfrak{l})U(\mathfrak{l}_+)}\mathbb{C}v_{\Omega,\psi},$$
where $Z(\mathfrak{l})U(\mathfrak{l}_+)=Z(\mathfrak{l})\ot U(\mathfrak{l}_+)$ acts by the characters $\Omega$ and $\psi$ on the $1$-dimensional space spanned by the vector $v_{\Omega,\psi}$.  As the Whittaker character $\psi$ is nonzero on each root space $\fg_\ga$ with $\ga\in\Delta_\psi$, it is easy to show that $Y_{\Omega,\psi}$ is a simple left $U(\mathfrak{l})$-module.  It remains irreducible when restricted to the (semisimple) derived subalgebra $\mathfrak{s}=[\mathfrak{l},\mathfrak{l}]$ of $\mathfrak{l}$.  We regard $Y_{\Omega,\psi}$ as a left $U(\mathfrak{p})$-module via the inflation map $\mathfrak{p}\rightarrow\mathfrak{p}/\mathfrak{m}=\mathfrak{l}$.  Inducing to $U(\mathfrak{g})$ gives the Whittaker module
$$M_{\Omega,\psi}=U(\mathfrak{g})\ot_{U(\mathfrak{p})}Y_{\Omega,\psi}.$$
\begin{theorem}{\em \cite[Theorem 2.9]{mcdowell85}}  The module $M_{\Omega,\psi}$ has a unique simple quotient $L_{\Omega,\psi}$, and every simple Whittaker module of type $\psi$ is isomorphic to a module $L_{\Omega,\psi}$ for some character $\Omega:\ Z(\mathfrak{l})\rightarrow\mathbb{C}$.\qed

  \end{theorem}

\section{Whittaker vectors}

We maintain the notation of the previous section.  Recall that the centre $\mathfrak{z}$ of the Lie algebra $\mathfrak{l}=\mathfrak{l}_\psi$ is contained in the Cartan subalgebra $\mathfrak{h}$ of $\mathfrak{g}$.  Define a partial order $\preceq$ on its algebraic dual $\mathfrak{z}^*$ by setting $\eta\preceq\mu$ if $\mu-\eta$ is the restriction of an element $\lambda\in\hbox{Span}_{\mathbb{Z}_+}(\Delta\setminus\Delta_\psi)$ to $\mathfrak{z}$.  We now state our main result.
  \begin{theorem}\label{mainthm}
Let $V$ be a Whittaker module for $\fg$.  Then $V$ is simple if and only if $\hbox{Wh}\,V$ is $1$-dimensional.
    \end{theorem}

  \noindent
  \proof Let $\psi:\ \eta\rightarrow\mathbb{C}$ be the Whittaker character of $V$.  The module $V$ has an obvious filtration $V=\bigcup_{k=0}^\infty V^{(k)}$, where
  $$V^{(k)}=\{v\in V:\ x_0\bullet x_1\bullet\cdots\bullet x_k\bullet v=0\hbox{\ for all\ }x_0,x_1,\ldots ,x_k\in\mathfrak{n}\},$$
  and $x\bullet v$ is defined as $xv-\psi(x)v$ for all $x\in\mathfrak{n}$ and $v\in V$. 

  Let $v$ be a nonzero element of $V$.  If $v\notin\hbox{Wh}\,V$, then $x\bullet v\neq0$ for some $x\in\mathfrak{n}$.  But $x\bullet v\in V^{(k-1)}$ whenever $v\in V^{(k)}$, so $V^{(0)}\cap U(\fg)v\neq 0$.   If $\hbox{Wh}\,V$ is $1$-dimensional and spanned by a vector $w$, then by definition, $V^{(0)}=\hbox{Wh}\,V=\mathbb{C}w$ and $V=U(\mathfrak{g})w$.  In this case, we see that $\mathbb{C}w= V^{(0)}\subseteq U(\fg)v$.  Since $w$ generates $V$, the module $V$ is thus simple when $\hbox{Wh}\,V$ is $1$-dimensional.

  Conversely, by McDowell's classification, any simple Whittaker module is of the form $L_{\Omega,\psi}$ for some Whittaker character $\psi:\ \mathfrak{n}\rightarrow \mathbb{C}$ and central character $\Omega\in Z(\mathfrak{l}_\psi)^*$.  The centre $\mathfrak{z}$ of the Lie algebra $\mathfrak{l}=\mathfrak{l}_\psi$ acts semisimply on $Y_{\Omega,\psi}$.  Since $\mathfrak{z}\subseteq\mathfrak{h}$, it also acts semisimply on the induced module $M_{\Omega,\psi}=U(\mathfrak{g})\ot_{U(\mathfrak{p})}Y_{\Omega,\psi}$ and on its simple quotient $L_{\Omega,\psi}$.  That is,
  $$M_{\Omega,\psi}=\bigoplus_{\mu\in\mathfrak{z}^*}M_{\Omega,\psi}^\mu\hbox{\ and\ }L_{\Omega,\psi}=\bigoplus_{\mu\in\mathfrak{z}^*}L_{\Omega,\psi}^\mu,$$
  where $V^\mu=\{v\in V\ :\ hv=\mu(h)v\hbox{\ for all\ }h\in\mathfrak{z}\}$ for any $\fg$-module $V$.  Since $\mathfrak{z}$ acts on $Y_{\Omega,\psi}$ by the restriction $\overline{\Omega}$ of the character $\Omega$ to $\mathfrak{z}$, it acts on
  $$M_{\Omega,\psi}=U(\overline{\mf{m}})\ot_\mathbb{C}U(\mf{p})\ot_{U(\mf{p})}Y_{\Omega,\psi}=U(\overline{\mf{m}})\ot_\mathbb{C} Y_{\Omega,\psi}$$
  by weights bounded above by $\overline{\Omega}$.  Similarly, 
\begin{equation}\label{star}
  L_{\Omega,\psi}=\bigoplus_{\mu\preceq\overline{\Omega}}L_{\Omega,\psi}^\mu.
\end{equation}

Let $v$ be a Whittaker vector in $L_{\Omega,\psi}$, with decomposition $v=\sum_{\mu\preceq\overline{\Omega}}v_\mu$ relative to (\ref{star}).  Let $x\in\mf{n}=\mf{l}_+\oplus\mf{m}$ be a root vector in a root space $\fg_\ga$.  If $x\in\mf{l}_+$, then $\mf{z}$ commutes with $x$ and $xL_{\Omega,\psi}^\mu\subseteq L_{\Omega,\psi}^\mu$.  In particular, $\sum_\mu\psi(x)v_\mu=xv=\sum_\mu xv_\mu$, so by comparing graded components, $xv_\mu=\psi(x)v_\mu$ for all $\mu$.  If $x\in\mf{m}$, then
$$zx=x(z+\ga(z)),$$
for all $z\in\mf{z}$, so
$$xL_{\Omega,\psi}^\mu\subseteq L_{\Omega,\psi}^{\mu+\overline{\ga}},$$
where $\overline{\ga}$ is the restriction of $\ga$ to $\mf{z}$.  But $\psi$ vanishes on $\mf{m}$ by construction, so
$$0=xv=\sum_{\mu\preceq\overline{\Omega}}xv_\mu.$$
Each $xv_\mu$ belongs to a distinct graded component, so $xv_\mu=0$ for all $\mu$.  The subalgebras $\mf{l}_+$ and $\mf{m}$ are spanned by root vectors, so $xv_\mu=\psi(x)v_\mu$ for all $x\in\mf{n}$ and $\mu\preceq\overline{\Omega}$.  That is, $v_\mu\in\hbox{Wh}\,V$ for all $\mu$.

Let $\mu\prec\overline{\Omega}$.  Then by degree considerations, the submodule
$$U(\fg)v_\mu=U(\overline{\mf{m}}\oplus\mf{l})v_\mu\subseteq\sum_{\lambda\preceq\mu}L_{\Omega,\psi}^\lambda$$
is proper in $L_{\Omega,\psi}$.  As $L_{\Omega,\psi}$ is simple, we see that $v_\mu=0$ and $v=v_{\overline{\Omega}}$.  Therefore, $v$ is a Whittaker vector in the simple $\mf{s}$-Whittaker module $Y_{\Omega,\psi}=L_{\Omega,\psi}^{\overline{\Omega}}$ of type $\overline{\psi}$, where $\overline{\psi}$ is the restriction of $\psi$ to the positive part $\mf{l}_+$ of $\mf{s}=[\mf{l},\mf{l}]$.  As $Y_{\Omega,\psi}$ is non-singular, Kostant's criterion \cite[Theorem 3.6.1]{kostant78} applies, and its space $W$ of Whittaker vectors is $1$-dimensional.  The space $\hbox{Wh}\,L_{\Omega,\psi}\subseteq W$ is thus also $1$-dimensional.\qed

\bigskip

Let $\mathcal{N}$ be the category of $\fg$-modules which are locally finite with respect to the action of $\mf{n}$.  An immediate consequence of Theorem \ref{mainthm} is that, for any module in $\mathcal{N}$, the dimension of its space of Whittaker vectors is bounded by its composition length.

\begin{corollary}
Let $M$ be a $\fg$-module in the category $\mathcal{N}$.  Then $\dim\hbox{Wh}\,M$ is bounded by the length of $M$.
  \end{corollary}

\noindent
\proof Without loss of generality, we may assume that $M$ is of finite length and has a composition series
$$0=M_0\subset M_1\subset\cdots\subset M_\ell=M.$$
We induct on the length $\ell$ of $M$, and as the base case $\ell=0$ is trivial, we may assume that $\ell>0$.

For any $i>0$ and nonzero $m\in M_i/M_{i-1}$, the Lie subalgebra $\mf{n}$ acts on the finite dimensional space $U(\mf{n})m$.  By Lie's theorem, $U(\mf{n})m\subseteq M_i/M_{i-1}$ contains a nonzero Whittaker vector $n$, and by simplicity $M_i/M_{i-1}=U(\fg)n$ is thus a Whittaker module.  The space $\hbox{Wh}(M_i/M_{i-1})$ is $1$-dimensional by Theorem \ref{mainthm}.

If $v,w\in M$ are nonzero Whittaker vectors, then their images $v+M_{\ell-1}$ and $w+M_{\ell-1}$ lie in the $1$-dimensional space $\hbox{Wh}(M_\ell/M_{\ell-1})$. The vector $v$ is thus of the form $\lambda w + u$ for some $\lambda\in\mathbb{C}$ and $u\in M_{\ell-1}$.  As $u=v-\lambda w$ is clearly also a Whittaker vector and $\dim\hbox{Wh}(M_{\ell-1})\leq \ell-1$ by induction, it follows that $\dim\hbox{Wh}\,M\leq 1+\dim\hbox{Wh}(M_{\ell-1})\leq\ell.$\qed

\bigskip

\begin{remark}{\em 
  Any Whittaker module $V$ of type $\psi$ that admits a central character $\chi:\ Z(\fg)\rightarrow \mathbb{C}$ is a quotient of a universal module
  $$\mathcal{Y}_{\chi,\psi}=U(\fg)\ot_{Z(\fg)U(\mf{n})}\mathbb{C}w_{\chi,\psi},$$
  where $Z(\fg)U(\mf{n})=Z(\fg)\ot U(\mf{n})$ acts on the $1$-dimensional space $\mathbb{C}w_{\chi,\psi}$ by $\chi\ot\psi$.  Such modules are always simple if $\psi$ is non-singular, but need not be simple for singular $\psi$.  By the above corollary, the dimension of $\hbox{Wh}\,V$ is bounded by the length of $\mathcal{Y}_{\chi,\psi}$.  Such lengths have been computed by Mili\v{c}i\'c and Soergel \cite{milicic-soergel97} using Kazhdan-Lusztig theory.
  }\end{remark}

\bigskip

\noindent
Karthik Dulam, Department of Mathematics, Indian Institute of Science, Bengaluru, India 560012, karthik.06.dulam@gmail.com

\medskip

\noindent
Hrishikesh Ghate, Department of Mathematics, Indian Institute of Science Education and Research, Bhopal, India 462066, hrishikesh19@iiserb.ac.in

\medskip

\noindent
Michael Lau, D\'epartement de math\'ematiques et de statistique, Universit\'e Laval, Qu\'ebec, Canada G1V0A6, Michael.Lau@mat.ulaval.ca

\medskip

\noindent
Suyash Pathak, Department of Mathematics and Statistics, Indian Institute of Technology, Kanpur, India 208016, skpat042001@gmail.com


\begin{thebibliography}{}

%\bibitem[AL]{fdreps} J.~Auger and M.~Lau, {Extensions of modules for twisted current algebras}, in preparation.

%\bibitem[BBT]{Babelon-Bernard-Talon03}
%O.~Babelon, D.~Bernard, and M.~Talon, {Introduction to Classical Integrable Systems}, Cambridge, 2003.  

  
%\bibitem[Ba]{batra}
%P. Batra, Representations of twisted multi-loop Lie algebras, J. Algebra {\bf 272} (2004), 404--416.


\bibitem{block80} R.~Block, The irreducible representations of the Lie algebra $\mathfrak{sl}(2)$ and of the Weyl algebra, Adv.~Math. {\bf 39} (1981), 69--110.

  
%\bibitem[Bog]{bogoyavlensky76} O.I.~Bogoyavlensky, On perturbations of the periodic Toda lattice, Commun.~Math.~Phys. {\bf 51} (1976), 201--209.

%\bibitem[Bo1]{bourbaki-algebre8} N. Bourbaki, \'El\'ements de math\'ematique: Alg\`ebre, Chapitre 8, Hermann, Paris, 1958.

%\bibitem[Bo2]{bourbaki-commalg5} N.~Bourbaki, \'El\'ements de math\'ematique: Alg\`ebre commutative, Chapitre 5, Hermann, Paris, 1964.

%\bibitem[Bo3]{bourbaki-Lie1} N. Bourbaki, \'El\'ements de math\'ematique: Groupes et alg\`ebres de Lie, Chapitre 1, Hermann, Paris, 1972.

%\bibitem[Bou]{bourbaki-Lie7-8} N.~Bourbaki, \'El\'ements de math\'ematique: Groupes et alg\`ebres de Lie, Chapitres 7 et 8, Hermann, Paris, 1975.

%\bibitem[BBL]{BBL}
%G.~Benkart, D.~Britten, and F.~Lemire, Modules with bounded weight multiplicities for simple {L}ie algebras, Math.~Z. {\bf 225} (1997), 333-353.

%\bibitem[BL]{BrLe}
%D.~Britten and F.~Lemire, A classification of simple {L}ie modules having a $1$-dimensional weight space, Trans.~Amer.~Math.~Soc. {\bf 299} (1987), 683--697.

%\bibitem[Ca]{Casati-Ortenzi06}
%  P.~Casati and G.~Ortenzi, New integrable hierarchies from vertex operator representations of polynomial Lie algebras, J.~Geom.~Phys. {\bf 56} (2006), 418--449.
  
%\bibitem[Ch]{chari86}
%V. Chari, Integrable representations of affine Lie-algebras, Invent. Math. {\bf 85} (1986), 317--335.


%\bibitem[CP]{Chari-Pressley94}
%V.~Chari and A.~Pressley, {A Guide to Quantum Groups}, Cambridge, 1994.  


  %\bibitem[CP]{cp88}
%V. Chari and A. Pressley, Integrable representations of twisted affine Lie algebras, J. Algebra {\bf 113} (1988), 438--464.


%\bibitem[CFK]{CFK}
%V.~Chari, G.~Fourier, and T.~Khandai, A categorical approach to Weyl modules, Transform. Groups {\bf 15} (2010), 517--549.
%\bibitem[DR]{dateroan00}
%E. Date and S.-s. Roan, The structure of quotients of the Onsager algebra by closed ideals, J. Phys. A {\bf 33} (2000), 3275--3296.
%\bibitem[EW]{erdmann-wildon}
%K. Erdmann and M. Wildon, Introduction to Lie Algebras, Springer Undergrad. Math. Ser., Springer, London, 2006.
%\bibitem[FL]{feigenloktev04}
%B. Feigin and S. Loktev, Multi-dimensional Weyl modules and symmetric functions, Comm. Math. Phys. {\bf 251} (2004), 427--445.

%\bibitem[DG]{DG}
%I.~Dimitrov and D.~Grantcharov, Classification of simple weight modules over affine Lie algebras, arXiv preprint 2009: 0910.0688v1.

%\bibitem[DMP]{DMP}
%I.~Dimitrov, O.~Mathieu, and I.~Penkov, On the structure of weight modules, Trans.~Amer.~Math.~Soc. {\bf 352} (2000), 2857--2869.

%\bibitem[D]{dixmier}
%  J.~Dixmier, {Alg\`ebres enveloppantes}, \'Editions Gauthier-Villars, Paris, 1974. 


%\bibitem[Fe]{fernando}
%S.~Fernando, Lie algebra modules with finite dimensional weight spaces, I, Trans.~Amer.~Math.~Soc. {\bf 322} (1990), 757--781.

%\bibitem{GP} P.~Gille and A.~Pianzola, {\it Galois cohomology and forms of algebras over Laurent polynomial rings}, Math. Ann. {\bf 338} (2007), 497--543.

%\bibitem[GLZ]{GLZ}
%X.~Guo, R.~Lu, and K.~Zhao, Simple Harish-Chandra modules, intermediate series modules, and Verma modules over the loop-Virasoro algebra, Forum Math. {\bf 23} (2011, 1029--1052.

%\bibitem[H]{hartwig07}
%B. Hartwig, The tetrahedron algebra and its finite-dimensional irreducible modules, Linear Algebra Appl. {\bf 422} (2007), 219--235.

%\bibitem[Ka]{Ka85}
%C. Kassel, K\"ahler differentials and coverings of complex simple Lie algebras extended over a commutative algebra, J.~Pure Appl.~Algebra {\bf 34} (1984), 265--275.

%\bibitem{KnusOjanguren} M.-A.~Knus and M.~Ojanguren, {Th\'eorie de la descente et alg\`ebres d'Azumaya}, Lecture Notes in Math., vol. 389, Springer, Berlin, 1974.

\bibitem{kostant78} B.~Kostant, On Whittaker vectors and representation theory, Invent.~Math. {\bf 48} (1978), 101--184.
  
\bibitem{kostant79} B.~Kostant, The solution to a generalized Toda lattice and representation theory, Adv.~Math. {\bf 34} (1979), 195--338.

%\bibitem[L1]{multiloop} M.~Lau, Representations of multiloop algebras, Pacific J. Math. {\bf 245} (2010), 167--184.

%\bibitem[L2]{twcurr} M.~Lau, Representations of twisted current algebras, J.~Pure Appl.~Algebra {\bf 218} (2014), 2149--2163.

%\bibitem[LP]{repforms} M.~Lau and A.~Pianzola, Maximal ideals and representations of twisted forms of algebras, Algebra Number Theory {\bf 7} (2013), 431--448.

%\bibitem{losev10} I.~Losev, 

  
%\bibitem[Ma1]{mathieu92} O.~Mathieu, Classification of Harish-Chandra modules over the Virasoro Lie algebra, Invent.~Math. {\bf 107} (1992), 225-234.

%\bibitem[Ma]{Ma00} O.~Mathieu, Classification of irreducible weight modules, Ann.~Inst.~Fourier, Grenoble {\bf 50} (2000), 537--592.

%\bibitem{matsumura} H.~Matsumura, {Commutative Ring Theory}, Cambridge, 1989.

\bibitem{mcdowell85} E.~McDowell, On modules induced from Whittaker modules, J.~Algebra {\bf 96} (1985), 161--177.
  
\bibitem{milicic-soergel97} D.~Mili\v{c}i\'c and W.~Soergel, The composition series of modules induced from Whittaker modules, Comment.~Math.~Helv. {\bf 72} (1997), 503--520.
  
%\bibitem{milne} J.S.~Milne, {\'Etale Cohomology}, Princeton, 1980.

%\bibitem[M]{molev} A.I.~Molev, Casimir elements and Sugawara operators for Takiff algebras, J.~Math.~Phys. {\bf 62} (2021), no.~1, paper no.~011701. 

%\bibitem[NSS]{NSS}
%E.~Neher, A.~Savage, and P.~Senesi, Irreducible finite-dimensional representations of equivariant map algebras, Trans. Amer. Math. Soc. {\bf 364} (2012), 2619--2646.

%\bibitem{PiPrSu} A.~Pianzola,  D.~Prelat, and J.~Sun, {\it Descent constructions for central extensions of infinite dimensional Lie algebras}, Manuscripta Math.  {\bf 122}  (2007), 137--148.

\bibitem{premet02} A.~Premet, Special transverse slices and their enveloping algebras, Adv.~Math. {\bf 170} (2002), 1--55.
  
%\bibitem[Ra]{rao93}
%S.E. Rao, On representations of loop algebras, Comm. Algebra {\bf 21} (1993), 2131--2153.

%\bibitem[R2]{rao01}
%S.E. Rao, Classification of irreducible integrable modules for multi-loop algebras with finite-dimensional weight spaces, J. Algebra {\bf 246} (2001), 215--225.

%\bibitem[Sa]{savage}
%A.~Savage, Classification of irreducible quasifinite modules over map Virasoro algebras, Forum Math. {\bf 17} (2012), 547--570.

%\bibitem{senesi}
%P.~Senesi, {\em Finite-dimensional representation theory of loop algebras: a survey}, Quantum Affine Algebras, Extended Affine Algebras, and their Applications (Banff, Canada, 2008), Contemp. Math., vol. 506, Amer. Math. Soc., 2010, pp. 263--283.

\bibitem{skryabin02} S.~Skryabin, A category equivalence, Appendix to \cite{premet02}.
  
  
%\bibitem[T]{toda67} M.~Toda, Wave propagation in anharmonic lattices, J.~Phys.~Soc.~Japan {\bf 23} (1967), 501--506.


%\bibitem[V]{vinberg} E.B.~Vinberg, On certain commutative subalgebras of a universal enveloping algebra (Russian), Izv.~Akad.~Nauk SSSR Ser.~Mat. {\bf 54} (1990), 3--25; translation in Math.~USSR-Izv. {\bf 36} (1991), 1--22.

  
\end{thebibliography}
\end{document}